\RequirePackage{ifpdf}
\ifpdf 
\documentclass[pdftex]{sigma}
\else
\documentclass{sigma}
\fi

\begin{document}

\renewcommand{\PaperNumber}{077}

\FirstPageHeading

\ShortArticleName{Global Stability of Dynamic Systems of High
Order}

\ArticleName{Global Stability of Dynamic Systems of High Order}

\Author{Mohammed BENALILI and Azzedine LANSARI}

\AuthorNameForHeading{M. Benalili and A. Lansari}

\Address{Department of Mathematics, B.P.~119, Faculty of
Sciences,\\ University Abou-bekr BelKa\"{\i}d, Tlemcen, Algeria}
\Email{\href{mailto:m_benalili@mail.univ-tlemcen.dz}{m\_benalili@mail.univ-tlemcen.dz},
\href{mailto:a_lansari@mail.univ-tlemcen.dz}{a\_lansari@mail.univ-tlemcen.dz}}

\ArticleDates{Received December 18, 2006, in f\/inal form June 04,
2007; Published online July 15, 2007}

\Abstract{This paper deals with global asymptotic stability of
prolongations of f\/lows induced by specif\/ic vector f\/ields and
their prolongations. The method used is based on various estimates
of the f\/lows.}


\Classification{37C10; 34D23}

\section{Introduction}

Global stability of dynamic systems is a vast domain in ordinary
dif\/ferential equations and it is one of its main topics. Many
works have been done in this context, we list some  of them:
\cite{3,4,5,6,7,8}. However, little is known in the stability of
high order (see \cite{10} and \cite{2}). In this paper, we are
concerned with the global asymptotic stability of prolongations of
f\/lows generated by some specif\/ic vector f\/ields and their
perturbations. The method used is based on various estimates of
the f\/lows and their prolongations. To justify the study of the
dynamic of prolongations of f\/lows, we consider the Lie algebra
$\chi ({\mathbb R}^{n})$ of vector f\/ields on ${\mathbb R}^{n}$
endowed with the weak topology, which is the topology of the
uniform convergence of vector f\/ields and all their derivatives
on a compact sets. The Lie bracket is a fundamental operation not
only in dif\/ferential geometry but in many f\/ields of
mathematics, such  as dynamic and control theory. The
invertibility of this latter is of many uses i.e.\ given any
vector f\/ields $X$, $Z$ f\/ind a vector f\/ield~$Y$ such that
$\left[ X,Y\right] =Z$. In the case of vector f\/ields $X$
def\/ined in a neighborhood of a~point~$a$ with $X(a)\neq 0$ we
have a positive answer: since in this case the vector f\/ield $X$
is locally of the form $\frac{\partial }{\partial x_{1}}$ and the
solution is given by
\begin{gather*}
Y(x_{1},\dots ,x_{n})=\int_{-r}^{x_{1}}Z(t,x_{2},\dots ,x_{n})dt,
\end{gather*}
where $\left\Vert x\right\Vert =\max\limits_{1\leq i\leq
n}\left\vert x_{i}\right\vert <r$. In the case of singular vector
f\/ields, i.e.\ $X(a)=0$ little is known. Consider a singular
vector f\/ield $X$ def\/ined in a neighborhood $U$ of the origin
$0$ with $X(0)=0$ and let $\phi _{t}$ be the f\/low generated by
$X$. Suppose that $X$ is complete and consider a vector
f\/ield~$Y$ def\/ined on an open set $V$ $\supset \phi _{t}(U)$ for all $%
t\in {\mathbb R}$. The transportation of a vector f\/ield $Y$
along the f\/low $\phi _{t}$ is def\/ined as
\begin{gather*}
(\phi _{t})_{\ast }Y(x)=\left( D\phi _{t}\cdot Y\right) \circ \phi
_{-t}(x)
\end{gather*}
and the derivative with respect to $t$ is given as follows
\begin{gather*}
\frac{d}{dt}(\phi _{t})_{\ast }Y=\left[ (\phi _{t})_{\ast }X,(\phi
_{t})_{\ast }Y\right] .
\end{gather*}
Put $Y_{t}=-\int_{0}^{t}(\phi _{s})_{\ast }Zds$, then
\begin{gather*}
\left[ X\text{,}Y_{t}\right] =-\frac{d}{dt}\big| _{t=0}(\phi
_{t})_{\ast }\int_{0}^{t}(\phi _{s})_{\ast }Zds
=-\int_{0}^{t}\frac{d}{ds}(\phi _{s})_{\ast }Zds=Z-(\phi
_{t})_{\ast }Z.
\end{gather*}
So if $(\phi _{t})_{\ast }Z$ converges to $0$ and the integral $%
Y=-\int_{0}^{+\infty }(\phi _{s})_{\ast }Zds$ is convergent in the
weak topology, then $Y$ is a solution of our equation.

As applications of the right invertibility of the bracket
operation on germs of vector f\/ields at a singular point we refer
the reader to the papers by the authors \cite{1,2} (see also
\cite{10}).

\section{Generalities}

First we recall some def\/initions on global asymptotic stability
as introduced in~\cite{9}. Let $\left\Vert \cdot\right\Vert $ be
the Euclidean norm on ${\mathbb R}^{n}$, $K\subset {\mathbb
R}^{n}$ is a compact set and $f$ any smooth function on ${\mathbb
R}^{n}$, we put
\begin{gather}
\left\Vert f\right\Vert _{r}^{K}=\sup_{x\in K}\max_{\left\vert
\alpha \right\vert \leq r}\left\Vert D^{\alpha }f(x)\right\Vert.
\label{1}
\end{gather}

\begin{definition}
A point $a\in {\mathbb R}^{n}$ is said \ globally asymptotically
stable (in brief $G.A.S.$) of the f\/low $\phi _{t}$ if

i) $a$ is an asymptotically stable (in brief $A.S.$) equilibrium
of the f\/low $\phi _{t}$;

ii) for any compact set $K\subset {\mathbb R}^{n}$ and any
$\varepsilon >0$ there exists $T_{K}>0$ such that $\ $for any
$t\geq T_{K}\ $we have $\Vert \phi _{t}\left( x\right) -a\Vert
\leq \epsilon $ for all $x\in K$.
\end{definition}

\begin{definition}
The point $a\in {\mathbb R}^{n}$ is said globally asymptotically
stable of order $r$ ($1\leq r\leq \infty $) for the f\/low $\phi
_{t}$ if

i) $a$ is a $G.A.S.$ point for the f\/low $\phi _{t}$;

ii) for any compact set $K\subset {\mathbb R}^{n}$ and
\begin{gather*}
\forall \; \epsilon >0 , \ \ \exists \; T_{K}>0 \ \ \text{such
that}\ \ \forall\; t\geq T_{K}\ \Rightarrow \  \Vert \phi
_{t}-aI\Vert _{r}^{K}\leq \epsilon ,
\end{gather*}
where $I$ denotes the identity map.
\end{definition}

A vector f\/ield $X$ will be called semi-complete if the
$X$-f\/low $\phi _{t}=\exp (tX)$ is def\/ined for all $t\geq 0$.

First we quote the following proposition which characterizes the
uniform asymptotic stability, for a proof see the book of
W.~Hahn~\cite{5}.

Let $(\phi )_{t}$ denote a f\/low def\/ined on ${\mathbb R}^{n}$.

\begin{proposition}
\label{prop1} The origin $0$ in ${\mathbb R}^{n}$ is $G.A.S.$
point for the flow $\phi _{t}$ if  for any ball $B(0,\rho )$,
centered at $0$ and of radius $\rho >0$, there exist $t_0\geq 0$
and functions $a$, $b$ such that
\begin{gather}
\Vert \phi _{t}(x)\Vert \leq a(\Vert x\Vert )b(t)   \label{2}
\end{gather}
with $a$ a continuous function on $B(0,\rho )$ monotonously
increasing such that $a(0)=0$ and $b$ is a~continuous function
defined for any $t\geq t_0$ monotonously decreasing such that
$\lim\limits_{t\rightarrow +\infty } b(t)=0$.
\end{proposition}

\section[Estimates of prolongations of flows]{Estimates of prolongations of f\/lows}

We start with some perturbations of linear vector f\/ields.

\subsection[Perturbation of linear vector field]{Perturbation of linear vector f\/ields}

Consider the following linear vector f\/ield
\begin{gather*}
X_{1}=\sum_{i=1}^{n}\alpha _{i}x_{i}\frac{\partial }{\partial
x_{i}},
\end{gather*}
where the coef\/f\/icients $\alpha _{i}$ $\in $ $\left[ a,b\right]
\subset {\mathbb R}$ and are not all $0$.

The $X_{1}$-f\/low, $\psi _{t}^{1}=\exp (tX_{1})$ \ is then
\begin{gather}
\psi _{t}^{1}(x)=xe^{\alpha t}=\left( x_{1}e^{\alpha _{1}t},\dots
,x_{n}e^{\alpha _{n}t}\right) \quad \forall\; t\in {\mathbb R}
\label{3}
\end{gather}
and its estimates are given by
\begin{gather}
\left\Vert x\right\Vert e^{at}\leq \Vert \psi _{t}^{1}(x)\Vert
\leq \left\Vert x\right\Vert e^{bt}.   \label{4}
\end{gather}

Consider now a perturbation of the vector f\/ield $X_{1}$ of the form $%
Y_{1}=X_{1}+Z_{1}$, where $Z_{1}$ is a~smooth vector f\/ield
globally Lipschitzian on ${\mathbb R}^{n}$. The explicit form of
the $Y_{1}$-f\/low is then
\begin{gather}
\psi _{t}^{1}(x)=xe^{At}+\int_{0}^{t}Z_{1}\left( \psi
_{s}^{1}(x)\right) ds,
 \label{5}
\end{gather}
where $A=\left(
\begin{array}{ccc}
\alpha _{1} & \cdots & 0 \\
\vdots & \ddots & \vdots \\
0 & \cdots & \alpha _{n}
\end{array}
\right) $.

\begin{lemma}
\label{lem1} If the perturbation $Z_{1}$ fulfills
\begin{gather}
\left\Vert Z_{1}(x)\right\Vert \leq c_{0}\quad \forall\; x\in
{\mathbb R}^{n}
 \label{6}
\end{gather}
then the vector field $Y_{1}$ is complete and the $Y_{1}$-flow
satisfies the estimates
\begin{gather*}
\left(\left\Vert x\right\Vert
-\frac{c_{0}}{a}\right)e^{bt}+\frac{c_{0}}{a}\leq \Vert \psi
_{t}^{1}(x)\Vert \leq \left(\left\Vert x\right\Vert +\frac{
c_{0}}{b}\right)e^{bt}-\frac{c_{0}}{b}.
\end{gather*}
\end{lemma}

\begin{proof}
Clearly the $Y_{1}$-f\/low $\psi _{t}^{1}$ is bounded for any $t\in \left[ 0,T%
\right] $ with $T<+\infty $ \ and any $x\in {\mathbb R}^{n}$. The
same is true if we replace $t$ by $-t$. Then $\psi _{t}^{1}$ is
complete.

Consider now the equation
\begin{gather}
\frac{1}{2}\frac{d}{dt}\left\Vert \psi _{t}^{1}(x)\right\Vert
^{2}=\left\langle \psi _{t}^{1}(x),\alpha \psi
_{t}^{1}(x)+Z_{1}\left( \psi _{t}^{1}(x)\right) \right\rangle .
\label{7}
\end{gather}
Letting $y=\Vert \psi _{t}^{1}(x)\Vert $, we deduce
\begin{gather*}
ay^{2}-c_0y\leq \frac{1}{2}\frac{d}{dt}y^{2}\leq by^{2}+c_0
y,\qquad y(0)=\left\Vert x\right\Vert
\end{gather*}
and by integrating we obtain
\begin{gather*}
\left(\left\Vert x\right\Vert
-\frac{c_0}{a}\right)e^{bt}+\frac{c_0}{a}\leq y\leq
\left(\left\Vert x\right\Vert
+\frac{c_0}{b}\right)e^{bt}-\frac{c_0}{b}.\tag*{\qed}
\end{gather*}
  \renewcommand{\qed}{}
\end{proof}

Let $B(0,1)$ be the open unit ball centered at the origin $0$.

\begin{lemma}
\label{lem2} If the perturbation $Z_{1}$ fulfills the estimates
\begin{gather}
\left\Vert Z_{1}(x)\right\Vert \leq c_0^{\prime }\left\Vert
x\right\Vert ^{1+m}\quad \forall \; x\in B\left( 0,1\right) \
\text{and any integer} \ m\geq
1, \nonumber \\
\left\Vert Z_{1}(x)\right\Vert \leq c_{0}^{\prime \prime
}\left\Vert x\right\Vert\quad \text{for every}  \ x\in {\mathbb
R}^{n}\setminus B\left(0,1\right),
 \label{8}
\end{gather}
then $Y_{1}$ is complete and the $Y_{1}$-flow fulfills the
following estimates for ant $t\geq 0$
\begin{gather}
\left\Vert x\right\Vert e^{a_0 t}\leq \left\Vert \psi
_{t}^{1}(x)\right\Vert \leq \left\Vert x\right\Vert e^{b_0t},\nonumber \\
\left\Vert x\right\Vert e^{-b_0 t}\leq \left\Vert \psi
_{-t}^{1}(x)\right\Vert \leq \left\Vert x\right\Vert e^{-a_0 t}
 \label{9}
\end{gather}
with $c_0=\max \left\{ c_0^{\prime }, c_0^{\prime \prime
}\right\}$, $a_0=a-c_0$ and $b_0=b+c_0$.
\end{lemma}

\begin{proof}
Taking account of the explicit form of the f\/low~(\ref{5}) and
the estimates (\ref{8}), we deduce that $Y_{1}$ is complete. If
$x\in B\left( 0,1\right) $ then $\left\Vert Z_{1}(x)\right\Vert
\leq c_0^{\prime }\left\Vert x\right\Vert ^{1+m}$ $\leq
c_0^{\prime }\left\Vert x\right\Vert $, letting $c_0=\max \left\{
c_0^{\prime },c_0^{\prime \prime }\right\} $ then $ \left\Vert
Z_{1}(x)\right\Vert \leq c_0\left\Vert x\right\Vert $ for any $
x\in {\mathbb R}^{n}$. If we put $y=\left\Vert \psi
_{t}^{1}(x)\right\Vert $ the equation~(\ref{7}) leads to
\begin{gather*}
(a-c_0)y\leq \frac{d}{dt}y\leq (b+c_0)y, \qquad y(0)=\left\Vert
x\right\Vert
\end{gather*}
and putting $b_0=b+c_0$, $a_0=a-c_0$, we deduce the following
estimates
\begin{gather*}
\left\Vert x\right\Vert e^{a_0t}\leq y\leq \left\Vert x\right\Vert
e^{b_0 t}\quad \text{for any} \ t\geq 0.
\end{gather*}
The same is also true in the on ${\mathbb R}^{n}\setminus
B\left(0,1\right) $.
\end{proof}

\begin{lemma}
\label{lem3} Suppose that all the coefficients $\alpha _{i}$ are negative, $%
a\leq \alpha _{i}\leq b<0$.

If the perturbation $Z_{1}$ fulfills the estimates
\begin{gather}
\left\Vert Z_{1}(x)\right\Vert \leq c_0\left\Vert x\right\Vert
^{1+m}\qquad \text{for any} \ x\in {\mathbb R}^{n}\ \text{and any
integer} \ m\geq 1,    \label{10}
\end{gather}
then the vector field $Y_{1}$ is semi-complete and the
$Y_{1}$-flow satisfies the estimates for any $t\geq 0$
\begin{gather}
\left\Vert x\right\Vert e^{at}\left( 1-\frac{c_0}{a}\left\Vert
x\right\Vert ^{m}(1-e^{amt})\right) ^{-\frac{1}{m}} \label{11}\\
\qquad{} \leq \left\Vert \psi _{t}^{1}(x)\right\Vert \leq \Vert
x\Vert e^{bt}\left( 1-\frac{c_0}{b}\left\Vert x\right\Vert
^{m}(1-e^{bmt})\right) ^{-\frac{1}{m}}.\nonumber
\end{gather}
\end{lemma}

\begin{proof}
By the relation~(\ref{5}) and the estimates~(\ref{10}), we deduce
that the vector f\/ield $Y_{1}$ is semi-complete. Letting
$y=\left\Vert \psi _{t}^{1}(x)\right\Vert $ and taking into
account the equation~(\ref{7}) and the estimates~(\ref{10}) we
deduce that
\begin{gather*}
ay-c_0y^{1+m}\leq \frac{d}{dt}y\leq by+c_0y^{1+m},\qquad
y(0)=\left\Vert x\right\Vert
\end{gather*}
and by integration we have
\begin{gather*}
\left\Vert x\right\Vert e^{at}\left( 1-\frac{c_0}{a}\left\Vert
x\right\Vert ^{m}(1-e^{amt})\right) ^{-\frac{1}{m}}\leq y\leq
\left\Vert x\right\Vert e^{bt}\left( 1-\frac{c_0}{b}\left\Vert
x\right\Vert ^{m}(1-e^{bmt})\right) ^{-\frac{1}{m}}.\tag*{\qed}
\end{gather*}
\renewcommand{\qed}{}
\end{proof}

\begin{example}
\label{exa1} Let the vector f\/ield
\begin{gather*}
X_{3}=\sum_{i=1}^{n}\left( \alpha _{i}x_{i}+\beta
_{i}x_{i}^{1+m_{i}}\right) \frac{\partial }{\partial x_{i}}
\end{gather*}
such that all the coef\/f\/icients fulf\/illing
\begin{gather*}
a\leq \alpha _{i}\leq b<0, \qquad a^{\prime }\leq \beta _{i}\leq
b^{\prime }\leq 0
\end{gather*}
and all the exponents $m_{i}$ are even positive integers with $%
0<m'_0\leq m_{i}\leq m_0$. The associated f\/low $\phi
_{t}^{3}=\exp (tX_{3})$ is the solution of the dynamic system
\begin{gather*}
\frac{d}{dt}\phi _{t}(x)=X_{3}\circ \phi _{t}(x),\qquad \phi _0
(x)=x
\end{gather*}
or in coordinates
\begin{gather*}
\frac{d}{dt}\left( \phi _{t}(x)\right) _{i}=\alpha _{i}\left( \phi
_{t}(x)\right) _{i}+\beta _{i}\left( \phi _{t}(x)\right)
_{i}^{1+m_{i}},\qquad \phi _0(x)=x.
\end{gather*}
This latter is a Bernoulli type equation and its solution is given
by
\begin{gather}
\left( \phi _{t}^{3}(x)\right) _{i}=x_{i}e^{\alpha _{i}t}\left(1+\frac{\beta _{i}%
}{\alpha _{i}}x_{i}^{m_{i}}\left( 1-e^{\alpha _{i}m_{i}t}\right) \right)^{\frac{-1%
}{m_{i}}}.    \label{12}
\end{gather}
The $X_{3}$-f\/low $\phi _{t}^{3}=\exp \left( tX_{3}\right) $ then
has the explicit form
\begin{gather*}
\phi _{t}^{3}(x)=xe^{\alpha t}\left(1+\frac{\beta }{\alpha
}x^{m}\left( 1-e^{\alpha mt}\right) \right)^{\frac{-1}{m}}
\end{gather*}
and the following estimates are true, $\forall\; t\geq 0$
\begin{gather}
\left\Vert x\right\Vert e^{at}\leq \left\Vert \phi
_{t}^{3}(x)\right\Vert \leq \left\Vert x\right\Vert
e^{bt}.\label{13}
\end{gather}

\end{example}

\subsection[Estimation of the $k^{\rm th}$ prolongation of the $Y_{1}$-flow]{Estimation of the $\boldsymbol{k^{\rm th}}$ prolongation of the $\boldsymbol{Y_{1}}$-f\/low}

Denote by $\eta _{1}^{1}(t,x,\nu )=D\psi _{t}^{1}(x)\nu $, where
$\nu \in {\mathbb R}^{n}$, the f\/irst derivative with respect to
$x$ of the $Y_{1}$-f\/low, solution of the dynamic system
\begin{gather*}
\frac{d}{dt}\eta _{1}^{1}(t,x,\nu )=\left(
D_{y}X_{1}+D_{y}Z_{1}\right) \eta _{1}^{1}(t,x,\nu ),\qquad \eta
_{1}^{1}(0,x,\nu )=\nu
\end{gather*}
with $y=\psi _{t}^{1}(x)$.

\begin{lemma}
\label{lem4} If the perturbation $Z_{1}$ fulfills the estimate
\begin{gather}
\left\Vert DZ_{1}(x)\right\Vert \leq c_{1}\quad \text{for any} \
x\in {\mathbb R}^{n},
  \label{14}
\end{gather}
then the derivative of the $Y_{1}$-flow is complete and has the
following estimates, for any $t\geq 0$
\begin{gather}
e^{a_{1}t}\leq \left\Vert D\psi _{t}^{1}(x)\right\Vert \leq
e^{b_{1}t},\qquad e^{-b_{1}t}\leq \left\Vert D\psi
_{-t}^{1}(x)\right\Vert \leq e^{-a_{1}t}
  \label{15}
\end{gather}
with $a_{1}=a-c_{1}$ and $b_{1}=b+c_{1}$ .
\end{lemma}

\begin{proof}
Consider as in previous lemmas the following equation
\begin{gather}
\frac{1}{2}\frac{d}{dt}\left\Vert \eta _{1}^{1}(t,x,\nu
)\right\Vert ^{2}=\left\langle \eta _{1}^{1}(t,x,\nu ),\left(
\alpha +DZ_{1}\right) \eta _{1}^{1}(t,x,\nu )\right\rangle
\label{16}
\end{gather}
and put $z=\left\Vert \eta _{1}^{1}(t,x,\nu )\right\Vert $, so
\begin{gather}
(a-c_{1})z^{2}\leq \frac{1}{2}\frac{d}{dt}z^{2}\leq
(b+c_{1})z^{2},\qquad z(0)=\left\Vert \nu \right\Vert
  \label{17}
\end{gather}
and then
\begin{gather*}
\left\Vert \nu \right\Vert e^{a_{1}t}\leq z\leq \left\Vert \nu
\right\Vert e^{b_{1}t}\quad \text{for any}\  \ t\geq 0\ \
\text{and} \ \ \nu \in {\mathbb R}^{n}. \tag*{\qed}
\end{gather*}\renewcommand{\qed}{}
\end{proof}

\begin{lemma}
\label{lem5} If the perturbation $Z_{1}$ fulfils the estimates
\begin{gather*}
\big\Vert D^{l}Z_{1}(x)\big\Vert \leq c_{l}^{\prime }\left\Vert
x\right\Vert ^{1-l+m}\  \text{for any} \  x\in B\left( 0,1\right) \ \text{and all integers} \ m\geq 1, \\
\big\Vert D^{l}Z_{1}(x)\big\Vert \leq c_{l}^{\prime \prime
}\left\Vert x\right\Vert ^{1-l}\quad \forall\;  x\in {\mathbb
R}^{n}\setminus B\left( 0,1\right)
\end{gather*}
with $l=0,1$, then the first derivative of the $Y_{1}$-flow is
complete and is estimated
by, for any $t\geq 0$%
\begin{gather}
e^{a_{1}t}\leq \Vert D\psi _{t}^{1}(x)\Vert \leq e^{b_{1}t},\qquad
e^{-b_{1}t}\leq \Vert D\psi _{-t}^{1}(x)\Vert \leq e^{-a_{1}t}
 \label{18}
\end{gather}
with $c_{l}=\max \left\{ c_{l}^{\prime },c_{l}^{\prime \prime
}\right\} $, $ a_{l}=a-c_{l}$ and $b_{l}=b+c_{l}$, $l=0,1$.
\end{lemma}

\begin{proof}
For any $x\in B\left( 0,1\right) $ we have $\Vert
D^{l}Z_{1}(x)\Vert \leq c_{l}^{\prime }\left\Vert x\right\Vert
^{1-l+m} $ $\leq c_{l}^{\prime }\left\Vert x\right\Vert ^{1-l}$
and letting $c_{l}=\max \left\{ c_{l}^{\prime },c_{l}^{\prime
\prime }\right\} $, we get for any $x\in {\mathbb R}^{n}$ $\Vert
D^{l}Z_{1}(x)\Vert \leq c_{l}\left\Vert x\right\Vert ^{1-l}$. By
the same arguments as in previous lemmas we get the
estimates~(\ref{18}).
\end{proof}

\begin{lemma}
\label{lem6} Suppose that all the coefficients $\alpha _{i}$ are negative, $%
a\leq \alpha _{i}\leq b<0$.

If the perturbation $Z_{1}$ fulfills the estimates
\begin{gather*}
\left\Vert Z_{1}(x)\right\Vert \leq c_{0}\left\Vert x\right\Vert
^{1+m},\quad \left\Vert DZ_{1}(x)\right\Vert \leq c_{1}\left\Vert
x\right\Vert ^{m} \quad \text{for all}  \ x\in {\mathbb R}^{n}  \
\text{and any integers}  \ m\geq 1.
\end{gather*}
Then the estimates of the first derivation of the $Y_{1}$-flow are
as follows, for any $t\geq 0$
\begin{gather*}
e^{at}\left( 1-\frac{c_{0}}{a}\left\Vert x\right\Vert
^{m}(1-e^{amt})\right) ^{-\frac{c_{1}}{mc_{0}}} \leq \left\Vert
D\psi _{t}^{1}(x)\right\Vert \leq e^{bt}\left(
1-\frac{c_{0}}{b}\left\Vert x\right\Vert ^{m}(1-e^{bmt})\right)
^{-\frac{c_{1}}{mc_{0}}}.
\end{gather*}
\end{lemma}

\begin{proof}
Letting $y=\left\Vert \psi _{t}^{1}(x)\right\Vert $ and
$z=\left\Vert \eta _{1}^{1}(t,x,\nu )\right\Vert $ in
equation~(\ref{16}), we get
\begin{gather*}
(a-c_{1}y^{m})z^{2}\leq \frac{1}{2}\frac{d}{dt}z^{2}\leq
(b+c_{1}y^{m})z^{2},\qquad z(0)=\left\Vert \nu \right\Vert
\end{gather*}
and taking into account the estimates given by the
relation~(\ref{11}), we obtain
\begin{gather*}
\left\Vert x\right\Vert ^{m}e^{mat}\left(
1-\frac{c_{0}}{a}\left\Vert x\right\Vert ^{m}(1-e^{amt})\right)
^{-1} \leq y^{m} \leq \left\Vert x\right\Vert ^{m}e^{mbt}\left(
1-\frac{c_0}{b}\left\Vert x\right\Vert ^{m}(1-e^{bmt})\right)
^{-1}
\end{gather*}
consequently
\begin{gather*}
\left\Vert \nu \right\Vert \exp \left(
at-c_{1}\int_{0}^{t}\frac{\left\Vert x\right\Vert
^{m}e^{mas}ds}{1-\frac{c_{0}}{a}\left\Vert x\right\Vert
^{m}(1-e^{ams})}\right) \\
\qquad {}\leq z
\leq \left\Vert \nu \right\Vert \exp \left( bt+c_{1}\int_{0}^{t}\frac{%
\left\Vert x\right\Vert ^{m}e^{mbs}ds}{1-\frac{c_0}{b}\left\Vert
x\right\Vert ^{m}(1-e^{bms})}\right)
\end{gather*}
which has the solution
\begin{gather*}
\left\Vert \nu \right\Vert e^{at}\left( 1-\frac{c_0}{a}\left\Vert
x\right\Vert ^{m}(1-e^{amt})\right) ^{-\frac{c_{1}}{mc_{0}}}\\
\qquad {}\leq z \leq \left\Vert \nu \right\Vert e^{bt}\left(
1-\frac{c_0}{b}\left\Vert x\right\Vert ^{m}(1-e^{bmt})\right)
^{-\frac{c_{1}}{mc_{0}}}\quad  \text{for} \ \ \nu \in {\mathbb
R}^{n}.\tag*{\qed}
\end{gather*}\renewcommand{\qed}{}
\end{proof}

\begin{example}
\label{exa2} We consider the same vector f\/ield as in
Example~\ref{exa1}. Denote by $\xi $$_{3}^{1}(t,x,\nu )=D\phi
_{t}^{3}(x)\nu $, $\forall\,\nu \in {\mathbb R}^{n}$, the f\/irst
derivation of the $X_{3}$-f\/low. In coordinates, we have for any
$i,j=1,\dots,n$,
\begin{gather*}
\left( \phi _{t}^{3}(x)\right) _{i}=x_{i}e^{\alpha _{i}t}\left(1+\frac{\beta _{i}%
}{\alpha _{i}}x_{i}^{m_{i}}\left( 1-e^{\alpha _{i}m_{i}t}\right) \right)^{\frac{-1%
}{m_{i}}}
\end{gather*}
so we deduce that
\begin{gather*}
\frac{\partial }{\partial x_{j}}\left( \phi _{t}^{3}(x)\right)
_{i}=e^{\alpha _{i}t}\left(1+\frac{\beta _{i}}{\alpha
_{i}}x_{i}^{m_{i}}\left( 1-e^{\alpha _{i}m_{i}t}\right)
\right)^{-1-\frac{1}{m_{i}}}\delta _{j}^{i}
\end{gather*}
and by the estimates~(\ref{13}) we get
\begin{gather*}
e^{at}\leq \left\Vert D\phi _{t}^{3}(x)\right\Vert \leq e^{bt}.
\end{gather*}
The second derivative is
\begin{gather*}
\frac{\partial ^{2}}{\partial x_{i}^{2}}\left( \phi
_{t}^{3}(x)\right) _{i}=-(1+m_{i})\frac{\beta _{i}}{\alpha
_{i}}x_{i}^{-1+m_{i}}e^{\alpha
_{i}t}\left( 1-e^{\alpha _{i}m_{i}t}\right) \left(1+\frac{\beta _{i}}{\alpha _{i}}%
x_{i}^{m_{i}}\left( 1-e^{\alpha _{i}m_{i}t}\right) \right)^{-2-\frac{1}{m_{i}}}%
.
\end{gather*}
Consequently, for $l=1,2$ and any $x\in B\left( 0,\rho \right) $
with $\rho
>0$ arbitrary f\/ixed, there are constants $M_{l}>0$ such that
\begin{gather*}
\big\Vert D^{l}\phi _{t}^{3}(x)\big\Vert \leq M_{l}e^{bt}\text{.}
\end{gather*}
\end{example}

\subsection[Perturbation of a nonlinear vector field]{Perturbation of a nonlinear vector f\/ield}

Consider the nonlinear vector f\/ield
\begin{gather*}
X_{2}=\sum_{i=1}^{n}\beta _{i}x_{i}^{1+m_{i}}\frac{\partial }{\partial x_{i}}%
\quad \text{with all} \ \ m_{i}>0 \ \ \text{and all} \ \ \beta
_{i}\leq 0.
\end{gather*}

The explicit form of the $X_{2}$-f\/low is then given by
\begin{gather}
\phi _{t}^{2}(x)=x(1-m\beta tx^{m})^{\frac{-1}{m}}    \label{19}
\end{gather}
for any $t\geq 0$ in the sense
\begin{gather*}
\left( \phi _{t}^{2}(x)\right) _{i}=x_{i}(1-m_{i}\beta _{i}tx_{i}^{m_{i}})^{%
\frac{-1}{m_{i}}},\qquad 1\leq i\leq n.
\end{gather*}

\begin{lemma}
\label{lem7} If the following assumptions are true

i) all the coefficients $\beta $$_{i}$ are non positive,
$-a^{\prime }\leq \beta _{i}\leq -b^{\prime }\leq 0$

ii) all the exponents $m_{i}$ are even positive integers; $0<m_{0}\leq $$%
m_{i}\leq $ $m_{0}^{\prime }.$

Then the vector field $X_{2}$ is semi-complete and the
$X_{2}$-flow satisfies the estimates
\begin{gather}
\left\Vert x\right\Vert \left( 1+b' m_0 t\left\Vert x\right\Vert
^{m_0}\right) ^{\frac{-1}{m_0}}\leq \Vert \phi _{t}^{2}(x)\Vert
\leq \left\Vert x\right\Vert \left( 1+a'm'_0 t\left\Vert
x\right\Vert ^{m_0^{\prime }}\right) ^{\frac{-1}{m'_0}}\ \
\text{for any} \ \ t\geq 0.    \label{20}
\end{gather}
\end{lemma}

\begin{proof}
Clearly the f\/low $\phi _{t}^{2}=\exp (tX_{2})$ given by
(\ref{19}) is semi-complete i.e.\ def\/ined for all $t\geq 0$.

Consider the equation
\begin{gather*}
\frac{1}{2}\frac{d}{dt}\Vert \phi _{t}^{2}(x)\Vert
^{2}=\left\langle \phi _{t}^{2}(x),\beta \left( \phi
_{t}^{2}(x)\right) ^{1+m}\right\rangle
\end{gather*}
and put $y=\phi _{t}^{2}(x)$, then
\begin{gather*}
b^{\prime }y^{2+m_{0}}\leq \frac{1}{2}\frac{d}{dt}y^{2}\leq
a^{\prime }y^{2+m_{0}^{\prime }},\qquad y(0)=\left\Vert
x\right\Vert
\end{gather*}
and we get the estimates given in~(\ref{20}).
\end{proof}

\subsection[Estimation of the $k^{\rm th}$ order derivation of the $X_{2}$-flow]{Estimation of the $\boldsymbol{k^{\rm th}}$ order derivation of the $\boldsymbol{X_{2}}$-f\/low}

Let $\xi $$_{2}^{1}(t,x,\nu )=D\phi _{t}^{2}(x)\nu $, $\forall \;
\nu \in {\mathbb R}^{n}$ be the f\/irst derivation of the
$X_{2}$-f\/low.

By formula~(\ref{19}), we get in coordinates
\begin{gather*}
\frac{\partial }{\partial x_{j}}\left( \phi _{t}^{2}(x)\right)
_{i}=(1-m_{i}\beta _{i}tx_{i}^{m_{i}})^{-1-\frac{1}{m_{i}}}\delta _{i}^{j}%
\quad \text{with}\quad \delta _{i}^{j}=\left\{
\begin{array}{c}
1 \ \ \text{if}\ \ i=j, \\
0 \ \ \text{if}\ \ i\neq j,
\end{array}
\right.
\end{gather*}
where $i,j=1,\dots ,n$.

Consequently
\begin{gather}
\left( 1+b' mt\left\Vert x\right\Vert ^{m_0}\right) ^{-1-\frac{1}{%
m_0}}\leq \Vert D\phi _{t}^{2}(x)\Vert \leq
\big( 1+a'm'_0 t\left\Vert x\right\Vert ^{m'_0}\big) ^{-1-%
\frac{1}{m'_0}}.  \label{21}
\end{gather}
To get the estimates of the second derivative, we put
\begin{gather*}
w_{i}=1-m_{i}\beta _{i}tx_{i}^{m_{i}},
\end{gather*}
so
\begin{gather*}
\frac{d}{dx_{i}}w_{i}=m_{i}(w_{i}-1)x_{i}^{-1} \qquad
\text{and}\qquad \frac{\partial }{\partial x_{i}}\left( \phi
_{t}^{2}(x)\right) _{i}=w_{i}{}^{-1-\frac{1}{m_{i}}}.
\end{gather*}
Consequently
\begin{gather*}
\frac{\partial ^{2}}{\partial x_{i}^{2}}\left( \phi
_{t}^{2}(x)\right)
_{i}=(1+m_{i})x_{i}^{-1}w_{i}{}^{-\frac{1}{m_{i}}}\left(
w_{i}^{-2}-w_{i}^{-1}\right)
=x_{i}^{-1}w_{i}{}^{-\frac{1}{m_{i}}}\left(
\frac{a_{1}^{2}}{w_{i}}+\frac{a_{2}^{2}}{w_{i}^{2}}\right),
\end{gather*}
where $a_{1}^{2}$ and $a_{2}^{2}$ are real constants. Let $\rho
>0$ be any arbitrary and f\/ixed real number, then for any $x\in
B(0,\rho )$ and any
 $t\geq t_{0}>0$ and $l=1,2$ there is $M_{l}>0$ such that
\begin{gather*}
\big\Vert D^{l}\phi _{t}^{2}(x)\big\Vert \leq M_{l}t^{-1-\frac{1}{%
m'_0}}.
\end{gather*}
Suppose that for $l=1,\dots,k-1$, with f\/ixed $k$, there exist constants $%
a_{j}^{l}$ and $M_{l}>0$ such that
\begin{gather*}
\frac{\partial ^{l}}{\partial x_{i}^{l}}\left( \phi
_{t}^{2}(x)\right)_i
=x_{i}^{1-l}w_{i}{}^{-\frac{1}{m_{i}}}\sum_{j=1}^{l}\frac{a_{j}^{l}}{w_{i}^{j}},
\end{gather*}
where $a_{j}^{l}$ are real constants and
\begin{gather*}
\big\Vert D^{l}\phi _{t}^{2}(x)\big\Vert \leq M_{l}t^{-1-\frac{1}{%
m'_0}}\quad \forall \; t>0.
\end{gather*}
For the estimates of the $k^{\rm th}$ derivative, we compute
\begin{gather*}
\frac{\partial ^{k}}{\partial x_{i}^{k}}\left( \phi
_{t}^{2}(x)\right)
_{i}=x_{i}^{1-k}w_{i}{}^{-\frac{1}{m_{i}}}\sum_{j=1}^{k}\frac{a_{j}^{k}}{%
w_{i}^{j}},
\\
\frac{\partial ^{k}}{\partial x_{i}^{k}}\left( \phi _{t}^{2}(x)\right) _{i}=%
\frac{d}{dx_{i}}x_{i}^{2-k}w_{i}{}^{-\frac{1}{m_{i}}}\sum_{j=1}^{k-1}\frac{%
a_{j}^{k-1}}{w_{i}^{j}}
\\
\quad{}=x_{i}^{1-k}w_{i}{}^{-\frac{1}{m_{i}}}\sum_{j=1}^{k-1}\left( \frac{%
a_{j}^{k-1}}{w_{i}^{j}}(1-k-jm_{i})+\frac{a_{j}^{k-1}}{w_{i}^{j+1}}%
(1+jm_{i})\right)
=x_{i}^{1-k}w_{i}{}^{-\frac{1}{m_{i}}}\sum_{j=1}^{k}\frac{a_{j}^{k}}{%
w_{i}^{j}},
\end{gather*}
where $a_{j}^{k}$ are real constants.

So we resume

\begin{proposition}
\label{prop2} Suppose that

i) all the coefficients satisfy $\beta $$_{i}\leq 0$, $-a' \leq
\beta _{i}\leq -b'$,

ii) the exponents $m_{i}$ are even natural numbers such that
$0<m_{0}\leq m_{i}\leq m'_{0}$.

Let $\rho >0$ be any arbitrary fixed real number. For any $x\in
B(0,\rho )$, for any $t\geq t_{0}>0$ and $\forall\; k\geq 1$ there
exist a constant $M_{k}>0$ such that
\begin{gather}
\big\Vert D^{k}\phi _{t}^{2}(x)\big\Vert \leq M_{k}t^{-1-\frac{1}{%
m_{0}^{\prime }}}.  \label{22}
\end{gather}
\end{proposition}

\subsection[Estimates of the $Y_{2}$-flow]{Estimates of the $\boldsymbol{Y_{2}}$-f\/low}

Let
\begin{gather*}
Y_{2}=\sum_{i=1}^{n}\big( \beta _{i}x_{i}^{1+m_{i}}+Z_{2i}(x)\big) \frac{%
\partial }{\partial x_{i}}
\end{gather*}
the perturbation of the nonlinear vector f\/ield $X_{2}$ and
denote by $\psi _{t}^{2}=\exp (tY_{2})$ the solution of the
dynamic system
\begin{gather*}
\frac{d}{dt}\psi _{t}^{2}(x)=Y_{2}\circ \psi _{t}^{2}(x),\qquad
\psi _0^{2}(x)=x.
\end{gather*}
In coordinates we have, $i=1,\dots ,n$,
\begin{gather*}
\frac{\partial }{\partial t}\psi _{2,i}(t,x)=\beta _{i}\psi
_{2,i}^{1+m_{i}}(t,x)+Z_{2i}\left( \psi _{t}^{2}(x)\right),\qquad
\psi _{2,i}(0,x)=x_{i}.
\end{gather*}
Putting
\begin{gather*}
y_{i}(t)=\psi _{2,i}^{-m_{i}}(t,x)
\end{gather*}
and
\begin{gather*}
\psi _{t}^{2}(x)=y^{\frac{-1}{m}}(t)=\left( y_{1}^{\frac{-1}{m_{1}}%
}(t),\dots ,y_{n}^{\frac{-1}{m_{n}}}(t)\right)
\end{gather*}
we get
\begin{gather*}
y_{i}^{\prime }(t)=-m_{i}\psi _{2,i}^{-1-m_{i}}(t,x)\frac{\partial }{%
\partial t}\psi _{2,i}(t,x)\text{.}
\end{gather*}
The Cauchy problem reads as
\begin{gather*}
y_{i}^{\prime }(t)=-m_{i}\beta _{i}-m_{i}\left( y_{i}(t)\right) ^{1+\frac{1}{%
m_{i}}}Z_{2i}(y^{\frac{-1}{m}}(t)),\qquad y_{i}(0)=x_{i}^{-m_{i}}
\end{gather*}
and has the following solution
\begin{gather*}
y_{i}(t)=x_{i}^{-m_{i}}-m_{i}\beta _{i}t-m_{i}\int_{0}^{t}y_{i}(s)^{^{1+%
\frac{1}{m_{i}}}}Z_{2i}(y^{\frac{-1}{m}}(s)ds,
\end{gather*}
i.e.
\begin{gather*}
\psi _{2,i}(t,x)=x_{i}\left(1-m_{i}\beta
_{i}tx_{i}^{m_{i}}-m_{i}x_{i}^{m_{i}}\int_{0}^{t}\psi
_{i}(s,x)^{-1-m_{i}}Z_{2i}(\psi
_{s}^{2}(x))ds\right)^{-\frac{1}{m_{i}}},
\end{gather*}
so we have the explicit form of the $Y_{2}$-f\/low
\begin{gather}
\psi _{t}^{2}(x)=x\left(1-m\beta tx^{m}-mx^{m}\int_{0}^{t}\psi
_{s}^{2}(x)^{-1-m}Z_{2}(\psi _{s}^{2}(x))ds\right)^{-\frac{1}{m}}.
\label{23}
\end{gather}
Now we will estimate the $Y_{2}$-f\/low.

\begin{lemma}
\label{lem8} Suppose that

$i)$ all the coefficients satisfy $\beta $$_{i}\leq 0$,
$-a^{\prime }\leq \beta _{i}\leq -b^{\prime }$;

$ii)$ the exponents $m_{i}$ are even natural numbers with
$0<m_{0}\leq m_{i}\leq m_{0}^{\prime }$;

$iii)$
\begin{gather*}
\left\Vert Z_{2i}(x)\right\Vert \leq c_{0}^{\prime }\left\vert
x_{i}\right\vert ^{2+m_{i}}\quad \text{if} \ \ x\in B\left( 0,1\right), \\
\left\Vert Z_{2i}(x)\right\Vert \leq c_{0}^{\prime \prime
}\left\vert x_{i}\right\vert ^{1+m_{i}}\quad \text{if} \ \ x\in
{\mathbb R}^{n}\setminus B\left( 0,1\right)
\end{gather*}
with $c_0=\max \left\{ c_0^{\prime },c_0^{\prime \prime }\right\} $, $%
b_0=b^{\prime }-c_0>0$, $a_0=a^{\prime }+c_0$.

Then

$1)$ the vector field $Y_{2}$ is semi-complete;

$2)$ the $Y_{2}$-flow has the estimates
\begin{gather}
\left\Vert x\right\Vert \left( 1+a_0 m_0t\left\Vert x\right\Vert
^{m_0}\right) ^{\frac{-1}{m_0}}\leq \Vert \psi _{t}^{2}(x)\Vert
\leq \left\Vert x\right\Vert \big(
1+b_0m_0^{\prime }t\left\Vert x\right\Vert ^{m_0^{\prime }}\big) ^{%
\frac{-1}{m_0^{\prime }}};   \label{24}
\end{gather}

$3)$ let $\rho >0$ and $t_0>0$ be fixed, then for any $x\in
B\left( 0,\rho \right) $ and any $t\geq t_0>0$ there is a~constant
$M_0>0$ such that
\begin{gather}
\Vert \psi _{t}^{2}(x)\Vert \leq M_0\left\Vert x\right\Vert t^{-%
\frac{1}{m_{0}^{\prime }}}.  \label{25}
\end{gather}
\end{lemma}

\begin{proof}
Let $x\in B\left( 0,1\right) $, by assumption we have $\left\Vert
Z_{2i}(x)\right\Vert \leq c_0^{\prime }\left\vert x_{i}\right\vert
^{2+m_{i}}\leq c_0^{\prime }\left\vert x_{i}\right\vert ^{1+m_{i}}$, put $%
c_0=\max \left\{ c_0^{\prime },c_0^{\prime \prime }\right\} $ then
for any $x\in {\mathbb R}^{n}$ we deduce $\left\Vert
Z_{2i}(x)\right\Vert \leq c_0\left\vert x_{i}\right\vert
^{1+m_{i}}$. Now taking account of the relation~(\ref{23}) we
deduce that for any $t\in \left[ 0,T\right] $
\begin{gather*}
\left\Vert \psi _{t}^{2}(x)\right\Vert \leq \left\Vert
x\right\Vert \left(
1+mt\left\Vert x\right\Vert ^{m}(b^{\prime }-c_0)\right) ^{-\frac{1}{m}%
}\leq \left\Vert x\right\Vert
\end{gather*}
hence the vector $Y_{2}$ is semi-complete, i.e.\ def\/ined for all
$t\geq 0$.

Consider the equation
\begin{gather*}
\frac{1}{2}\frac{d}{dt}\Vert \left( \psi _{t}^{2}(x)\right)
_{i}\Vert ^{2}=\left\langle \left( \psi _{t}^{2}(x)\right)
_{i},\beta _{i}\left( \psi _{t}^{2}(x)\right)
_{i}^{1+m_{i}}+Z_{2i}\left( \psi _{t}^{2}(x)\right) \right\rangle
\end{gather*}
we get $y_{i}=\Vert \left( \psi _{t}^{2}(x)\right) _{i}\Vert $ and
$y_{i}(0)=\left\vert x_{i}\right\vert $, so we deduce
\begin{gather*}
\frac{1}{2}\frac{d}{dt}y_{i}^{2}\leq (\beta
_{i}+c_0)y_{i}^{2+m_{i}}\leq -(b^{\prime }-c_0)y_{i}^{2+m_{i}}
\end{gather*}
and
\begin{gather*}
\frac{1}{2}\frac{d}{dt}y_{i}^{2}\geq (\beta
_{i}-c_0)y_{i}^{2+m_{i}}\geq -(a^{\prime
}+c_0)y_{i}^{2+m_{i}}\text{.}
\end{gather*}
We put $b_0=b^{\prime }-c_0$ and $a_0=a^{\prime }+c_0$, the
solutions are estimated as
\begin{gather}
(\left\vert x_{i}\right\vert
^{-m_{i}}+a_{0}m_{i}t)^{-\frac{1}{m_{i}}}\leq \Vert \left( \psi
_{t}^{2}(x)\right) _{i}\Vert \leq (\left\vert x_{i}\right\vert
^{-m_{i}}+b_0m_{i}t)^{-\frac{1}{m_{i}}}. \label{26}
\end{gather}

Hence, we have the estimate~(\ref{25}).
\end{proof}

Now, we estimate the f\/irst derivation of the $Y_{2}$-f\/low. Let $\eta $$%
_{2}^{1}(t,x,\nu )=D\psi _{t}^{2}(x)\nu $, $\forall \; \nu \in
{\mathbb R}^{n}$ the solution of the dynamic system
\begin{gather*}
\frac{d}{dt}\eta _{2}^{1}(t,x,\nu )=\left(
D_{y}X_{2}+D_{y}Z_{2}\right) \eta _{2}^{1}(t,x,\nu ), \qquad \eta
_{2}^{1}(0,x,\nu )=\nu
\end{gather*}
with $y=\psi _{t}^{2}(x)$.

\begin{lemma}
\label{lem9} Suppose that

$i)$ the coefficients are such that $\beta $$_{i}\leq 0$,
$-a^{\prime }\leq \beta _{i}\leq -b^{\prime }$;

$ii)$ the coefficients $m_{i}$ are even natural numbers,
$0<m_0\leq m_{i}\leq m_0^{\prime }$;

$iii)$
\begin{gather*}
\Vert D^{l}Z_{2i}(x)\Vert \leq c_{l}^{\prime }\left\vert
x_{i}\right\vert ^{2-l+m_{i}} \ \ \text{if} \ x\in B\left( 0,1\right), \\
\Vert D^{l}Z_{2i}(x)\Vert \leq c_{l}^{\prime \prime }\left\vert
x_{i}\right\vert ^{1-l+m_{i}} \ \ \text{if} \ x\in {\mathbb
R}^{n}\setminus B\left( 0,1\right)
\end{gather*}
with $l=0,1$;

$iv)$
\begin{gather*}
a_0=a^{\prime }+c_0, \qquad b_0=b^{\prime }-c_0>0
\end{gather*}
and
\begin{gather*}
a_{1}=a^{\prime }(1+m_0)+c_{1}, \qquad b_{1}=b^{\prime
}(1+m_0)-c_{1}>0
\end{gather*}
 with $c_{l}=\max \left\{ c_{l}^{\prime },c_{l}^{\prime \prime
}\right\}$.

Then the first derivation of the $Y_{2}$-flow has the following
estimates, for any $t>0$
\begin{gather}
\big(1+b_0m_0t\left\Vert x\right\Vert ^{m_{0}}\big)^{-\frac{a_{1}}{b_{0}m_{0}}%
}\leq \Vert D\psi _{t}^{2}(x)\Vert \leq \big(1+a_0m_0^{\prime
}t\left\Vert x\right\Vert ^{m_0^{\prime }}\big)^{-\frac{b_{1}}{%
a_0m_0^{\prime }}}.  \label{27}
\end{gather}
Let $\rho >0$ be arbitrary and fixed for any $x\in B(0,\rho )$,
and any $t\geq t_0>0$ there is a constant $M_{1}>0$ such that
\begin{gather}
\Vert D\psi _{t}^{2}(x)\Vert \leq M_{1}t^{-\frac{b_{1}}{%
a_{0}m_{0}^{\prime }}}.  \label{28}
\end{gather}
\end{lemma}

\begin{proof}
Let $x\in B\left( 0,1\right) $, for $l=0,1$ we have
\begin{gather*}
\Vert D^{l}Z_{2i}(x)\Vert \leq c_{l}^{\prime }\left\vert
x_{i}\right\vert ^{2-l+m_{i}}\leq c_{l}^{\prime }\left\vert
x_{i}\right\vert ^{1-l+m_{i}}.
\end{gather*}
Let $c_{l}=\max \left\{ c_{l}^{\prime },c_{l}^{\prime \prime
}\right\} $ then for $x\in {\mathbb R}^{n}$ one has
\begin{gather*}
\Vert D^{l}Z_{2i}(x)\Vert \leq c_{l}\left\vert x_{i}\right\vert
^{1-l+m_{i}}.
\end{gather*}
Consider the equation
\begin{gather*}
\frac{1}{2}\frac{d}{dt}\Vert \eta _{2}^{1}(t,x,\nu )\Vert
^{2}=\left\langle \eta _{2}^{1}(t,x,\nu ),\left(
D_{y}X_{2}+D_{y}Z_{2}\right) \eta _{2}^{1}(t,x,\nu )\right\rangle
\end{gather*}
and put $z(t)=\Vert \eta _{2}^{1}(t,x,\nu )\Vert $ with $%
z(0)=\left\Vert \nu \right\Vert $, then
\begin{gather*}
\frac{1}{2}\frac{d}{dt}z^{2}\leq \sup_{i=1,\dots,n}\left(
((1+m_{i})\beta _{i}+c_{1})\Vert \left( \psi _{t}^{2}(x)\right)
_{i}\Vert ^{m_{i}}\right) z^{2} \leq z^{2}\sup_{i=1,\dots,n}\left(
-b_{1}\Vert \left( \psi _{t}^{2}(x)\right) _{i}\Vert
^{m_{i}}\right)
\end{gather*}
and
\begin{gather*}
\frac{1}{2}\frac{d}{dt}z^{2}\geq \inf_{i=1,\dots,n}\left(
((1+m_{i})\beta _{i}-c_{1})\Vert \left( \psi _{t}^{2}(x)\right)
_{i}\Vert ^{m_{i}}\right) z^{2} \geq z^{2}\inf_{i=1,\dots,n}\left(
-a_{1}\Vert \left( \psi _{t}^{2}(x)\right) _{i}\Vert
^{m_{i}}\right).
\end{gather*}
The solutions fulf\/ill the following estimates
\begin{gather*}
\left\Vert \nu \right\Vert \exp \inf_{i=1,\dots,n}\left(
-a_{1}\int_{0}^{t}\Vert \left( \psi _{s}^{2}(x)\right) _{i}\Vert
^{m_{i}}ds\right) \\
\qquad {}\leq z(t) \leq \left\Vert \nu \right\Vert \exp
\sup_{i=1,\dots,n}\left( -b_{1}\int_{0}^{t}\Vert \left( \psi
_{s}^{2}(x)\right) _{i}\Vert ^{m_{i}}ds\right)
\end{gather*}
with, by (\ref{26})
\begin{gather*}
\frac{\left\vert x_{i}\right\vert ^{m_{i}}}{1+a_0m_{i}t\left\vert
x_{i}\right\vert ^{m_{i}}}\leq \Vert \left( \psi
_{t}^{2}(x)\right)
_{i}\Vert ^{m_{i}}\leq \frac{\left\vert x_{i}\right\vert ^{m_{i}}}{%
1+b_0m_{i}t\left\vert x_{i}\right\vert ^{m_{i}}}.
\end{gather*}
So we deduce
\begin{gather*}
\left\Vert \nu \right\Vert \exp \inf_{i=1,\dots,n}\left( -a_{1}\int_{0}^{t}%
\frac{\left\vert x_{i}\right\vert
^{m_{i}}}{1+b_{0}m_{i}s\left\vert
x_{i}\right\vert ^{m_{i}}}ds\right) \\
\qquad {} \leq z(t)
\leq \left\Vert \nu \right\Vert \exp \sup_{i=1,\dots,n}\left( -b_{1}\int_{0}^{t}%
\frac{\left\vert x_{i}\right\vert ^{m_{i}}}{1+a_0m_{i}s\left\vert
x_{i}\right\vert ^{m_{i}}}ds\right).
\end{gather*}
Consequently the solutions satisfy
\begin{gather*}
\left\Vert \nu \right\Vert
\inf_{i=1,\dots,n}(1+b_0m_{i}t\left\vert x_{i}\right\vert
^{m_{i}})^{-\frac{a_{1}}{b_{0}m_{i}}}\leq z(t) \leq \left\Vert \nu
\right\Vert \sup_{i=1,\dots,n}(1+a_0m_{i}t\left\vert
x_{i}\right\vert ^{m_{i}})^{-\frac{b_{1}}{a_{0}m_{i}}}.
\end{gather*}
Then there are constants $m_0>0$ and $m_0^{\prime }>0$ such that
\begin{gather*}
\left\Vert \nu \right\Vert (1+b_0m_0t\left\Vert x\right\Vert ^{m_{0}})^{-%
\frac{a_{1}}{b_0m_0}}\leq \Vert D\psi _{t}^{2}(x)\nu \Vert
\\
\qquad{} \leq \left\Vert \nu \right\Vert \big(1+a_0m_0^{\prime
}t\left\Vert x\right\Vert ^{m_0^{\prime
}}\big)^{-\frac{b_{1}}{a_0m_0^{\prime }}} \quad \forall \; \nu \in
{\mathbb R}^{n} \ \text{and for any} \ t>0.
\end{gather*}
Hence, we have the estimate~(\ref{28}).
\end{proof}

\subsection[Perturbation of binomial vector fields]{Perturbation of binomial vector f\/ields}

Let
\begin{gather*}
Y_{3}=\sum_{i=1}^{n}\big( \alpha _{i}x_{i}+\beta
_{i}x_{i}^{1+m_{i}}+Z_{3i}(x)\big) \frac{\partial }{\partial
x_{i}}
\end{gather*}
with $a\leq \alpha _{i}\leq b<0$, $a^{\prime }\leq \beta _{i}\leq
b^{\prime }\leq 0$ and $0<m_0\leq m_{i}\leq m_0^{\prime }$, be the
perturbation of the binomial vector f\/ield $X_{3}$ and let $\psi
_{t}^{3}=\exp (tY_{3})$ be the $Y_{3}$-f\/low which is the
solution of the dynamic system
\begin{gather*}
\frac{d}{dt}\psi _{t}(x)=Y_{3}\circ \psi _{t}(x),\qquad \psi
_0(x)=x
\end{gather*}
and in coordinates, we get
\begin{gather*}
\frac{\partial }{\partial t}\psi _{3,i}(t,x)=\alpha _{i}\psi
_{3,i}(t,x)+\beta _{i}\psi _{3,i}^{1+m_{i}}(t,x)+Z_{3,i}\left(
\psi _{t}^{3}(x)\right),\qquad \psi _{i}(0,x)=x_{i}
\end{gather*}
which is a Bernoulli type equation and by the same method as in
the proof of previous lemmas and with putting
\begin{gather*}
y_{i}(t)=\psi _{3,i}^{-m_{i}}(t,x)
\end{gather*}
and
\begin{gather*}
\psi _{t}^{3}(x)=y^{\frac{-1}{m}}(t)=\left( y_{1}^{\frac{-1}{m_{1}}%
}(t),\dots ,y_{n}^{\frac{-1}{m_{n}}}(t)\right)
\end{gather*}
we get the solution
\begin{gather*}
\psi _{3,i}(t,x)=
x_{i}e^{^{\alpha _{i}t}}\Bigg( 1+\frac{\beta _{i}}{\alpha _{i}}%
x_{i}^{m_{i}}(1-e^{^{\alpha _{i}m_{i}t}})\\
\phantom{\psi _{3,i}(t,x)=}{}-m_{i}x_{i}^{m_{i}}\int_{0}^{t}%
\left[ \psi _{3,i}(s,x)\right] ^{-1-m_{i}}Z_{3,i}(\psi
_{s}^{3}(x))e^{^{\alpha _{i}m_{i}s}}ds\Bigg) ^{\frac{-1}{m_{i}}}
\end{gather*}
and the implicit form of the $Y_{3}$-f\/low reads as
\begin{gather}
\psi _{t}^{3}(x)=xe^{\alpha t}\left( 1+\frac{\beta }{\alpha }%
x^{m}(1-e^{\alpha mt})-mx^{m}\int_{0}^{t}\left[ \psi
_{s}^{3}(x)\right] ^{-1-m}Z_{3}(\psi _{s}^{3}(x))e^{\alpha
ms}ds\right) ^{-\frac{1}{m}}. \label{29}
\end{gather}

\subsection[Estimation of the $Y_{3}$-flow]{Estimation of the $\boldsymbol{Y_{3}}$-f\/low}

By the same arguments as in the previous, we get the following
estimates of the $Y_{3}$-f\/low.

\begin{lemma}
\label{lem10} If the following assumptions are true

$i)$ all the coefficients $\alpha $$_{i}$ are negative, $-a\leq
\alpha _{i}\leq -b<0$;

$ii)$ all the coefficients $\beta $$_{i}$ are non positive,
$-a^{\prime }\leq \beta _{i}\leq -b^{\prime }$;

$iii)$ the exponents $m_{i}$ are even natural numbers with $0<m_{0}\leq $$%
m_{i}\leq $ $m_{0}^{\prime }$;

$iv)$
\begin{gather*}
\left\Vert Z_{3i}(x)\right\Vert \leq c_0^{\prime }\left\vert
x_{i}\right\vert ^{2+m_{i}} \ \ \text{if} \ x\in B\left( 0,1\right), \\
\left\Vert Z_{3i}(x)\right\Vert \leq c_0^{\prime \prime
}\left\vert x_{i}\right\vert ^{1+m_{i}} \ \ \text{if}\ x\in
{\mathbb R}^{n}\setminus B\left(0,1\right)
\end{gather*}
with $c_0=\max \left\{ c_0^{\prime },c_0^{\prime \prime }\right\}
$, $b_0=b^{\prime }-c_0>0$, $a_0=a^{\prime }+c_0$.

Then

$1)$ there exist constants $m>0$ and $m^{\prime }>0$ such that the $Y_{3}-$%
flow has the estimates, $\forall\; t\geq 0$%
\begin{gather*}
\left\Vert x\right\Vert e^{-at}\left( 1+\frac{a_{0}}{a}\left\Vert
x\right\Vert ^{m}(1-e^{-amt})\right) ^{-\frac{1}{m}}\\
\qquad{}\leq \Vert \psi _{t}^{3}(x)\Vert \leq \left\Vert
x\right\Vert e^{-bt}\left( 1+\frac{b_{0}}{b}\left\Vert
x\right\Vert ^{m^{\prime }}(1-e^{-bm^{\prime }t})\right) ^{-\frac{1}{%
m^{\prime }}};
\end{gather*}

$2)$ for any $t>0$ there are positive constants $c_{1}$ and
$c_{2}$ such that
\begin{gather}
c_{1}\left\Vert x\right\Vert e^{-at}\leq \Vert \psi
_{t}^{3}(x)\Vert \leq c_{2}\left\Vert x\right\Vert e^{-bt};
\label{30}
\end{gather}

$3)$ the vector field $Y_{3}$ is semi-complete.
\end{lemma}

By similar calculations as in previous lemmas, we get the
following estimates to the f\/irst derivative of the
$Y_{3}$-f\/low.

\begin{lemma}
\label{lem11} Suppose that

$i)$ all the coefficients $\alpha $$_{i}$ are negative, $-a\leq
\alpha _{i}\leq -b<0$;

$ii)$ all the coefficients $\beta $$_{i}$are non positive,
$-a^{\prime }\leq \beta _{i}\leq -b^{\prime }$;

$iii)$ the exponents $m_{i}$ are even natural numbers such that
$0<m_{0}\leq m_{i}\leq m_{0}^{\prime }$;

$iv)${\samepage
\begin{gather*}
\Vert D^{l}Z_{3i}(x)\Vert \leq c_{l}^{\prime }\left\vert
x_{i}\right\vert ^{2-l+m_{i}} \ \ \text{if} \ x\in B\left( 0,1\right), \\
\Vert D^{l}Z_{3i}(x)\Vert \leq c_{l}^{\prime \prime }\left\vert
x_{i}\right\vert ^{1-l+m_{i}} \ \ \text{if} \ x\in {\mathbb
R}^{n}\setminus B\left( 0,1\right)
\end{gather*}
 with $l=0,1$;}

$v)$
\begin{gather*}
a_0=a^{\prime }+c_0,\qquad b_0=b^{\prime }-c_0>0
\end{gather*}
and
\begin{gather*}
a_{1}=a^{\prime }(1+m_0)+c_{1}, \qquad b_{1}=b^{\prime
}(1+m_0)-c_{1}>0
\end{gather*}
with $c_{l}=\max \left\{ c_{l}^{\prime },c_{l}^{\prime \prime
}\right\} $.

Then there exist constants $m>0$ and $m^{\prime }>0$ such that for any $%
t\geq 0$
\begin{gather*}
e^{-at}\left( 1+\frac{b_0}{b}\left\Vert x\right\Vert
^{m}(1-e^{-bmt})\right) ^{-\frac{a_{1}}{b_0m}}\\
\qquad{}\leq \Vert D\psi _{t}^{3}(x)\Vert \leq e^{-bt}\left(
1+\frac{a_0}{a}\left\Vert x\right\Vert ^{m^{\prime
}}(1-e^{-am^{\prime }t})\right) ^{-\frac{b_{1}}{a_0m^{\prime }}}
\end{gather*}
and for any $t\geq 0$, there is a constant $M_{1}>0$ such that
\begin{gather}
\left\Vert D\psi _{t}^{3}(x)\right\Vert \leq M_{1}e^{-bt}.
\label{31}
\end{gather}
\end{lemma}

\section[Global stability of prolongations of flows]{Global stability of prolongations of f\/lows}

With notations of the previous sections, we will give global
stability of some f\/lows.

\subsection[Global stability of the $Y_{1}$-flow]{Global stability of the $\boldsymbol{Y_{1}}$-f\/low}

\begin{lemma}
\label{lem12} Let the vector fields
\begin{gather*}
Y_{1}=\sum_{i=1}^{n}\left( \alpha _{i}x_{i}+Z_{1i}(x)\right)
\frac{\partial }{\partial x_{i}}
\end{gather*}
with the following assumptions

$i)$ all the coefficients are negative, $-a\leq \alpha _{i}\leq
-b<0$;

$ii)$
\begin{gather*}
\left\Vert Z_{1}(x)\right\Vert \leq c_0^{\prime }\left\Vert
x\right\Vert ^{1+m}\quad \forall \; x\in B\left( 0,1\right) \
\text{and} \ \forall m\geq 1,
\\
\left\Vert Z_{1}(x)\right\Vert \leq c_0^{\prime \prime }\left\Vert
x\right\Vert \quad \forall \; x\in {\mathbb R}^{n}\setminus
B\left( 0,1\right);
\end{gather*}

$iii)$ $b_0=b-c_0>0$, where $c_0=\max \left\{ c_0^{\prime },
c_0^{\prime \prime }\right\} $.

Then the origin $0$ is a globally asymptotically stable equilibrium to the $%
Y_{1}$-flow $\psi _{t}^{1}$ on ${\mathbb R}^{n}$.
\end{lemma}

\begin{proof}
Let $\psi _{t}^{1}=\exp (tY_{1})$ be the $Y_{1}$-f\/low, then by
the assumptions and the estimates given by Lemma~\ref{lem2} we get
that
\begin{gather*}
\left\Vert \psi _{t}^{1}(x)\right\Vert \leq \left\Vert
x\right\Vert e^{-b_0t}\quad \forall\; t\geq 0 \ \text{and} \
\forall x\in {\mathbb R}^{n}
\end{gather*}
and by Proposition~\ref{prop1}, the origin $0$ is $G.A.S.$ for
$\psi _{t}^{1}$ on ${\mathbb R}^{n}$.
\end{proof}

\begin{example}
\label{exa3} We consider the vector f\/ield
\begin{gather*}
X_{3}=\sum_{i=1}^{n}\left( \alpha _{i}x_{i}+\beta
_{i}x_{i}^{1+m_{i}}\right) \frac{\partial }{\partial x_{i}}
\end{gather*}
of Example~1 with $a\leq \alpha _{i}\leq b<0$, $a^{\prime }\leq
\beta _{i}\leq b^{\prime }\leq 0$. The $X_{3}$-f\/low $\phi
_{t}^{3}=\exp \left( tX_{3}\right) $ is then given by
\begin{gather*}
\phi _{t}^{3}(x)=xe^{\alpha t}\left(1+\frac{\beta }{\alpha
}x^{m}\left( 1-e^{\alpha mt}\right) \right)^{\frac{-1}{m}}.
\end{gather*}
Let $\rho >0$ be arbitrary and f\/ixed real number. By the
estimates (\ref{13}), we have
for any $x\in B(0,\rho )$ and any $t\geq t_0\geq 0$%
\begin{gather*}
\left\Vert \phi _{t}^{3}(x)\right\Vert \leq \left\Vert x\right\Vert e^{-bt}%
.
\end{gather*}
By Proposition\ref{prop1} the origin $0$ is a $G.A.S.$ for the
f\/low $\phi _{t}^{3}$ on ${\mathbb R}^{n}$.
\end{example}

\subsection[Global stability of the first prolongation of the $Y_{1}$-flow]{Global stability of the f\/irst prolongation of the $\boldsymbol{Y_{1}}$-f\/low}

\begin{lemma}
\label{lem13}With the same assumptions as in Lemma~{\rm
\ref{lem12}} and the following conditions
\begin{gather*}
\left\Vert DZ_{1}(x)\right\Vert \leq c_{1}^{\prime }\left\Vert
x\right\Vert ^{m}\quad \forall \; x\in B\left( 0,1\right) \
\text{and} \ \forall\;
m\geq 1, \\
\left\Vert DZ_{1}(x)\right\Vert \leq c_{1}^{\prime \prime }\quad
\forall\; x\in {\mathbb R}^{n}\setminus B\left( 0,1\right)
\end{gather*}
with $b_{1}=b-c_{1}>0$ and $c_{1}=\max \left\{ c_{1}^{\prime
},c_{1}^{\prime \prime }\right\} $.

Then the origin $0$ is a globally asymptotically stable for the
first prolongation of the $Y_{1}$-flow $\psi _{t}^{1}$ on
${\mathbb R}^{n}$.
\end{lemma}

\begin{proof}
By the estimates (\ref{18}) and the hypothesis we deduce that
\begin{gather*}
\Vert D\psi _{t}^{1}(x)\nu \Vert \leq \left\Vert \nu \right\Vert
e^{-b_{1}t}\quad \forall \; t>0, \ \ \forall\; \nu \in {\mathbb
R}^{n}
\end{gather*}
and by Proposition~\ref{prop1},we obtain that the origin $0$ is a
$G.A.S.$ equilibrium on ${\mathbb R}^{n}$ for $\eta
_{1}^{1}(t,x,v)=D\psi _{t}^{1}(x)\nu $.
\end{proof}

\subsection[Global stability of the $k^{\rm th}$ prolongation of the $Y_{1}$-flow]{Global stability of the
 $\boldsymbol{k^{\rm th}}$ prolongation of the $\boldsymbol{Y_{1}}$-f\/low}

Suppose that

i) all the coef\/f\/icients are negative, $-a\leq \alpha _{i}\leq
-b<0$;

ii) for any $l=1,\dots,k-1$%
\begin{gather*}
\Vert D^{l}Z_{1}(x)\Vert \leq c_{l}^{\prime }\left\Vert
x\right\Vert ^{1-l+m}\quad \text{for any} \ x\in B\left(
0,1\right) \ \text{and for
any integer} \ m\geq l-1, \\
\Vert D^{l}Z_{1}(x)\Vert \leq c_{l}^{\prime \prime }\quad
\forall\; x\in {\mathbb R}^{n}\setminus B\left( 0,1\right),\\
a_0=a+c_0,\qquad
b_0=b-c_0>0,\\
a_{1}=a+c_{1}, \qquad b_{1}=b-c_{1}>0
\end{gather*}
with $c_{l}=\max \left\{ c_{l}^{\prime },c_{l}^{\prime \prime
}\right\}$, $b_{l}=c_{l}$ $\forall\; l\geq 2$.

Put $\eta _{1}^{l}(t,x,\nu ,\dots ,\nu )=D^{k}\psi _{t}^{1}(x)\nu
^{k}$, where $\nu \in {\mathbb R}^{n}$. Since by Lemmas~\ref{12}
and~\ref{13} the origin $0$ is an $G.A.S.$ equilibrium for $\eta
_{1}^{l}$, with $l=0,1,$ on ${\mathbb R}^{n}$, we suppose that
this property remains true for $l=0,1,\dots ,k-1$ with $k\geq 2$
i.e.\ for any $\rho >0$ and any $x\in B(0,\rho )$ there exist constants $%
M_{l}>0$ such that for any $t\geq t_0>0$%
\begin{gather*}
\Vert D^{l}\psi _{t}^{1}(x)\Vert \leq M_{l}e^{-b_{1}t}\text{.}
\end{gather*}
We will show that the origin $0$ is a $G.A.S.$ equilibrium for
$\eta _{1}^{k} $ on ${\mathbb R}^{n}$. $\eta _{1}^{k}(t,x,\nu
,\dots,\nu )=D^{k}\psi _{t}^{1}(x)\nu ^{k}$ is solution of the
dynamic system
\begin{gather*}
\frac{d}{dt}\eta _{1}^{k}=D_{y}Y_{1}\cdot \eta
_{1}^{k}+G_{1}^{k}(t,x,\nu ) ,\qquad \eta _{1}^{k}(0,x,\nu
,\dots,\nu )=\nu
\end{gather*}
with $y=\psi _{t}^{1}(x)$ and
\begin{gather*}
G_{1}^{k}(t,x,\nu )=\sum_{l=2}^{k}D_{y}^{l}Y_{1}(y)\sum_{\underset{i_{j}>0}{%
i_{1}+\dots +i_{l}=k}}\left( \prod_{j=1}^{l}D^{i_{j}}\psi
_{t}^{1}(x)\nu ^{i_{j}}\right)
\\
\phantom{G_{1}^{k}(t,x,\nu )}{}
=\sum_{l=2}^{k-1}D_{y}^{l}Z_{1}(y)\sum_{\underset{i_{j}>0}{i_{1}+\dots
+i_{l}=k} }\left( \prod_{j=1}^{l}D^{i_{j}}\psi _{t}^{1}(x)\nu
^{i_{j}}\right) +D_{y}^{k}Z_{1}(y)\left( D\psi _{t}^{1}(x)\nu
\right) ^{k}.
\end{gather*}
Consequently we get
\begin{gather*}
\eta _{1}^{k}(t,x,\nu ,\dots ,\nu )=D\psi _{t}^{1}(x)\nu
+\int_{0}^{t}D\Psi _{t-s}^{1}(\psi _{s}^{1}(x))G_{1}^{k}(s,x,\nu
)ds.
\end{gather*}
The integral is well def\/ined at $s=0$, since
\begin{gather*}
\lim_{s\rightarrow 0^{+}}D\psi _{t-s}^{1}(\psi _{s}^{1}(x))=D\psi _{t}^{1}(x)%
\text{ }
\end{gather*}
and there exist constants $A_{l}>0$ such that
\begin{gather*}
\lim_{s\rightarrow 0^{+}}G_{1}^{k}(s,x,\nu
)=\sum_{l=2}^{k}A_{l}D_{y}^{l}Z_{1}(y)\nu ^{k}\text{.}
\end{gather*}
We will show that it converges uniformly with respect to $x$ as
$t+\infty $. Put
\begin{gather*}
I_{k}=\int_{0}^{t}\Vert D\psi _{t-s}^{1}(\psi _{s}^{1}(x))\Vert
\Vert G_{1}^{k}(s,x,\nu )\Vert ds.
\end{gather*}
Since $\Vert D^{l}Z_{1}(x)\Vert \leq c_{l}$ $\forall\; l\geq 1$, $
\forall \; x\in {\mathbb R}^{n}$, there are constants $b_{l}>0$
such that $\forall\; y\in {\mathbb R}^{n}$, $\Vert
D_{y}^{l}Y_{1}(y)\Vert \leq b_{l}$ and by the assumption of
recurrence there exist constants $M_{l}>0$ such that
\begin{gather*}
\Vert D^{l}\psi _{t}^{1}(x)\Vert \leq M_{l}e^{-b_{1}t}\quad
\forall \;t\geq 0.
\end{gather*}
We deduce that there is a constant $C_{k}>0$ such that
\begin{gather*}
I_{k}\leq
\sum_{l=2}^{k}b_{l}M_{l}\int_{0}^{t}e^{-b_{1}(t-s+sl)}ds\leq
C_{k}e^{-b_{1}t}.
\end{gather*}
So for any $x\in {\mathbb R}^{n}$ one has
\begin{gather*}
\lim_{t\rightarrow +\infty }I_{k}\leq
\sum_{l=2}^{k}M_{l}b_{l}\left\Vert \nu
\right\Vert ^{l}\int_{0}^{+\infty }e^{-b_{1}s}ds=\frac{1}{b_{1}}%
\sum_{l=2}^{k}M_{l}b_{l}\left\Vert \nu \right\Vert ^{l}
\end{gather*}
and the integral $I_{k}$ is uniformly convergent with respect to
$x\in {\mathbb R}^{n}$ as $t\rightarrow $ $+\infty $. Consequently
\begin{gather*}
\lim_{t\rightarrow +\infty }\Vert \eta _{1}^{k}\Vert
=\lim_{t\rightarrow +\infty }\Vert D\psi _{t}^{1}(x)\nu \Vert
+\int_{0}^{+\infty }\lim_{t\rightarrow +\infty }\Vert D\psi
_{t-s}^{1}(\psi _{s}^{1}(x))\Vert \Vert G_{1}^{k}(s,x,\nu )\Vert
ds=0
\end{gather*}
and there is a constant $M_{k}^{\prime }>0$ such that
\begin{gather*}
\Vert \eta _{1}^{k}\Vert \leq \Vert D\psi _{t}^{1}(x)\nu \Vert
+\int_{0}^{+\infty }\Vert D\psi _{t-s}^{1}(\psi _{s}^{1}(x))\Vert
\Vert G_{1}^{k}(s,x,\nu )\Vert ds\leq M_{k}^{\prime }\Vert \nu
\Vert ^{k}e^{-b_{1}t}.
\end{gather*}
This show by Proposition~\ref{prop1} that the origin $0$ is a
$G.A.S.$ equilibrium to $\eta _{1}^{k}$ on ${\mathbb R}^{n}$. We
formulate our proving as follows

\begin{proposition}
\label{prop3} Let $k\geq 0$ be any integer. The origin $0$ is a
$G.A.S.$
equilibrium of order $k$ for  the $Y_{1}$-flow and there is a constant $%
M_{k}>0$ such that $\forall\; t>0$
\begin{gather}
\Vert D^{k}\psi _{t}^{1}(x)\Vert \leq M_{k}e^{-b_{1}t},\qquad
\Vert D^{k}\psi _{-t}^{1}(x)\Vert \leq M_{k}e^{a_{1}t}.
\label{32}
\end{gather}
\end{proposition}

\section[Global stability of a flow generated by nonlinear perturbed vector
fields]{Global stability of a f\/low generated\\ by nonlinear
perturbed vector f\/ields}

First we will start with monomial vector f\/ields.

\subsection[Global stability of the $X_{2}$-flow]{Global stability of the $\boldsymbol{X_{2}}$-f\/low}

Let
\[
X_{2}=\sum_{i=1}^{n}\beta _{i}x_{i}^{1+m_{i}}\frac{\partial
}{\partial x_{i}}
\]
with

(i) all the coef\/f\/icients $\beta $$_{i}\leq 0$ such that
$-a^{\prime }\leq \beta _{i}\leq -b^{\prime }$;

(ii) all the exponents $m_{i}$ are even natural integers with $0<m_0\leq $$%
m_{i}\leq $ $m_0^{\prime }$.

Let $\phi _{t}^{2}=\exp \left( tX_{2}\right) $ be the
$X_{2}$-f\/low. By the estimations (\ref{19}) we obtain
\begin{gather*}
\left\Vert \phi _{t}^{2}(x)\right\Vert \leq \left\Vert
x\right\Vert \big( 1+a^{\prime }m_0^{\prime }t\left\Vert
x\right\Vert ^{m_0^{\prime }}\big) ^{\frac{-1}{m_0^{\prime }}}.
\end{gather*}

Let $\rho >0$ be arbitrary f\/ixed, for any $x\in B(0,\rho )$ and
any $t\geq t_0>0$ there is a constant $M_0>0$ such that
\begin{gather*}
\Vert \phi _{t}^{2}(x)\Vert \leq M_0\left\Vert x\right\Vert t^{-%
\frac{1}{m_{0}^{\prime }}}.
\end{gather*}
By Proposition~\ref{prop1}, the origin is a globally
asymptotically stable equilibrium to the f\/low $\phi _{t}^{2}$
on~${\mathbb R}^{n}$.

Let $l=1,2,\dots $ any positive integer. By
Proposition~\ref{prop2}, we have: for any f\/ixed $\rho >0$, and
all $x\in B(0,\rho )$ and $t\geq t_{0}>0$, there exist constants
$M_{l}>0$ and $M_{l}^{\prime }>0$ such that
\begin{gather*}
\Vert D^{l}\phi _{t}^{2}(x)\Vert \leq M_{l}\text{ }t^{-1-\frac{1}{%
m_{0}^{\prime }}}\qquad \text{and}\qquad \Vert D^{l}\phi
_{0}^{2}(x)\Vert \leq M_{l}^{\prime }.
\end{gather*}

So the origin $0$ is a $G.A.S.$ equilibrium for $D^{l}\phi
_{t}^{2}(x)$ on $ {\mathbb R}^{n}$.

Resuming our proving, we get

\begin{proposition}
Let $k\geq 0$ be any integer. Under the above conditions $(i)$ and
$(ii)$, the origin $0$ is a $G.A.S.$ of order $k$ for the
$X_{2}$-flow on ${\mathbb R}^{n}$.
\end{proposition}

\subsection[Global stability of high order of the $Y_{2}$-flow]{Global stability of high order of the
$\boldsymbol{Y_{2}}$-f\/low}

Let
\[
Y_{2}=\sum_{i=1}^{n}\big( \beta
_{i}x_{i}^{1+m_{i}}+Z_{2,i}(x)\big) \frac{\partial }{\partial
x_{i}}
\] be a
smooth vector f\/ield on ${\mathbb R}^{n}$ such that

i) all the coef\/f\/icients $\beta $$_{i}\leq 0$ are non negative with $%
-a^{\prime }\leq \beta _{i}\leq -b^{\prime }$;

ii) $m_{i}$ are even natural numbers with $0<m_{0}\leq $$m_{i}\leq $ $%
m_{0}^{\prime }$;

iii) for $k=0,\dots,1+m_{i}$%
\begin{gather*}
\Vert D^{k}Z_{2i}(x)\Vert \leq c_{k}^{\prime }\left\vert
x_{i}\right\vert ^{2-k+m_{i}}\quad \text{if} \ x\in B\left( 0,1\right); \\
\Vert D^{k}Z_{2i}(x)\Vert \leq c_{k}^{\prime \prime }\left\vert
x_{i}\right\vert ^{1-k+m_{i}}\quad \text{if} \ x\in {\mathbb
R}^{n}\setminus B\left(0,1\right);
\end{gather*}

iv) for any $k\geq 2+m_{i}$
\begin{gather*}
\| D^kZ_{2i}(x)\| \leq c_k;
\end{gather*}

v)
\begin{gather*}
a_{0}=a^{\prime }+c_{0},\qquad a_{1}=a^{\prime }(1+m_{0})+c_{1},
\\
b_{0}=b^{\prime }-c_{0}>0, \qquad b_{1}=b^{\prime
}(1+m_{0})-c_{1}>a_{0}m_{0}^{\prime }
\end{gather*}
 with $c_{k}=\max \left\{ c_{k}^{\prime },c_{k}^{\prime \prime
}\right\}$.

\begin{remark}
If $x\in B\left( 0,1\right) $ then $\Vert D^{k}Z_{2i}(x)\Vert \leq
c_{k}^{\prime }\left\vert x_{i}\right\vert ^{2-k+m_{i}}\leq
c_{k}^{\prime }\left\vert x_{i}\right\vert ^{1-k+m_{i}}$. Putting $%
c_{l}=\max \left\{ c_{l}^{\prime },c_{l}^{\prime \prime }\right\}
$, we deduce that for any $x\in {\mathbb R}^{n}$ have $\Vert
D^{k}Z_{2i}(x)\Vert \leq c_{k}\left\vert x_{i}\right\vert
^{1-k+m_{i}}$.
\end{remark}

\subsubsection[Global stability of the $Y_{2}$-flow on ${\mathbb R}^{n}$]{Global stability of
the $\boldsymbol{Y_{2}}$-f\/low on $\boldsymbol{{\mathbb R}^{n}}$}

Let $\psi _{t}^{2}=\exp \left( tY_{2}\right) $ be the $Y_{2}$-f\/low and let $%
\rho >0$ be arbitrary and f\/ixed, so by the estimates (\ref{25})
for all $x\in B(0,\rho )$ and all $t\geq t_0>0$ there is a
constant $M_0>0$ such that
\begin{gather*}
\Vert \psi _{t}^{2}(x)\Vert \leq M_0\left\Vert x\right\Vert t^{-%
\frac{1}{m_0^{\prime }}}\text{.}
\end{gather*}
So by Proposition~\ref{prop1}, the origin $0$ is a $G.A.S.$
equilibrium for the $Y_{2}$-f\/low $\psi _{t}^{2}$ on ${\mathbb
R}^{n}$.

\subsubsection[Global stability of prolongation of the $Y_{2}$-flow on ${\mathbb R}^{n}$]{Global
stability of prolongation of the $\boldsymbol{Y_{2}}$-f\/low on
$\boldsymbol{{\mathbb R}^{n}}$}

We proceed by recurrence. Since it is already true for $k=0$, we
suppose that for any $l=1,\dots,k-1$, with $k\geq 2$, the origin
$0$ is a $G.A.S.$ to $ D^{l}\psi _{t}^{2}(x)$ on ${\mathbb R}^{n}$
that is to say for any f\/ixed $\rho >0$, all $x\in B(0,\rho )$
and all $t\geq t_{0}>0$ there are constants $M_{l}>0$ such that
\begin{gather*}
\Vert D^{l}\psi _{t}^{2}(x)\Vert \leq M_{l}t^{-\frac{b_{1}}{%
a_{0}m_{0}^{\prime }}}\qquad \text{and} \qquad \Vert D^{l}\psi
_{0}^{2}(x)\Vert \leq M_{l}^{\prime }.
\end{gather*}

We will show that $0$ is a $G.A.S.$ for $D^{k}\psi _{t}^{2}(x)$ on
${\mathbb R}^{n}$.

Put $\eta _{2}^{k}(t,x,\nu ,\dots,\nu )=D^{k}\psi _{t}^{2}(x)\nu ^{k}$ $%
\forall \; \nu \in {\mathbb R}^{n}$ which is solution of the
dynamic system
\begin{gather*}
\frac{d}{dt}\eta _{2}^{k}=D_{y}Y_{2}\cdot \eta
_{2}^{k}+G_{2}^{k}(t,x,\nu ),\qquad \eta _{2}^{k}(0,x,\nu
,\dots,\nu )=\nu
\end{gather*}
with $y=\psi _{t}^{2}(x)$ and
\begin{gather*}
G_{2}^{k}(t,x,\nu )=\sum_{l=2}^{k}D_{y}^{l}Y_{2}(y)\sum_{\underset{i_{j}>0}{%
i_{1}+\dots +i_{l}=k}}\left( \prod_{j=1}^{l}D^{i_{j}}\psi
_{t}^{2}(x)\nu ^{i_{j}}\right).
\end{gather*}
By the method of the resolvent, we deduce
\begin{gather*}
\eta _{2}^{k}(t,x,\nu ,\dots ,\nu )=D\psi _{t}^{2}(x)\nu
+\int_{0}^{t}D\psi _{t-s}^{2}(\psi _{s}^{2}(x))G_{2}^{k}(s,x,\nu
)ds .
\end{gather*}
Clearly the integral
\begin{gather*}
I_{k}^{1}=\int_{0}^{1}\Vert D\psi _{t-s}^{2}(\psi
_{s}^{2}(x))\Vert \Vert G_{2}^{k}(s,x,\nu )\Vert ds
\end{gather*}
is well def\/ined at $s=0$ and $s=t$, since
\begin{gather*}
\lim_{s\rightarrow 0^{+}}D\psi _{t-s}^{2}(\psi _{s}^{2}(x))=D\psi
_{t}^{2}(x).
\end{gather*}
By the recurrent assumption $D^{l}\psi _{0}^{2}(x)$ are bounded
and there exist constants $A_{l}>0$ such that
\begin{gather*}
\lim_{s\rightarrow 0^{+}}\Vert G_{2}^{k}(s,x,\nu )\Vert \leq
\sum_{l=2}^{k}A_{l}\Vert D_{x}^{l}Y_{2}(x)\nu ^{l}\Vert .
\end{gather*}
In the same way
\begin{gather*}
\lim_{s\rightarrow t^{-}}D\psi _{t-s}^{2}(\psi
_{s}^{2}(x))=\text{identity}.
\end{gather*}

Now, we have to show that
\begin{gather*}
I_{k}^{2}=\int_{1}^{t}\Vert D\psi _{t-s}^{2}(\psi
_{s}^{2}(x))\Vert \Vert G_{2}^{k}(s,x,\nu )\Vert ds
\end{gather*}
converges uniformly on any compact set $K\subset {\mathbb R}^{n}$
as $t\rightarrow 0$.

Let $x\in K$, by the relations~(\ref{26}) and~(\ref{28}) we get for all $%
t\geq 0$
\begin{gather*}
\left\Vert x\right\Vert \left( 1+a_0m_0t\left\Vert x\right\Vert
^{m_0}\right) ^{\frac{-1}{m_0}}\leq \Vert \psi _{t}^{2}(x)\Vert
\leq \left\Vert x\right\Vert \big(
1+b_0m_0^{\prime }t\left\Vert x\right\Vert ^{m_0^{\prime }}\big) ^{%
\frac{-1}{m_0^{\prime }}},
\\
(1+b_0m_0t\left\Vert x\right\Vert ^{m_0})^{-\frac{a_{1}}{b_0m_0}%
}\leq \Vert D\psi _{t}^{2}(x)\Vert \leq \big(1+a_0m_0^{\prime
}t\left\Vert x\right\Vert ^{m'_0}\big)^{-\frac{b_{1}}{%
a_0m_0^{\prime }}}.
\end{gather*}
So $\left\Vert y\right\Vert =\Vert \psi _{t}^{2}(x)\Vert \leq
\left\Vert x\right\Vert $ and $\left\Vert D\psi _{t-s}^{2}(\psi
_{s}^{2}(x))\right\Vert $ is bounded. Since for any $x\in {\mathbb R}^{n}$ and any $%
l=1,\dots,1+m_{i}$, $\Vert D^{l}Z_{2i}(x)\Vert \leq
c_{l}\left\vert x_{i}\right\vert ^{1-l+m_{i}}$ then
$D_{y}^{l}Y_{2}(y)$ are bounded. Now by the assumption of
recurrence there exist constants $M_{l}>0$ such that for any $t>0$
\begin{gather*}
\Vert D^{l}\psi _{t}^{2}(x)\Vert \leq M_{l}t^{-\frac{b_{1}}{%
a_0m_0^{\prime }}}
\end{gather*}
with $a_0m_0^{\prime }<b_{1}$ i.e.\ $\frac{b_{1}}{a_0m_0^{\prime
}}>1$, and we deduce the existence of constants $C_{l}>0$ such
that
\begin{gather*}
\lim_{t\rightarrow +\infty }I_{k}^{2}\leq
\sum_{l=2}^{k}C_{l}\int_{1}^{+\infty }s^{-\frac{lb_{1}}{a_0m_0^{\prime }}%
}ds\leq \sum_{l=2}^{k}C_{l}\left( \frac{lb_{1}}{a_0m_0^{{\prime }}}%
-1\right) ^{-1}.
\end{gather*}

The integral $I_{k}^{2}$ converges uniformly on any compact
$K\subset {\mathbb R}^{n}$ as $t\rightarrow $ $+\infty $.

Now since the integral is well def\/ined at $s=0$, then
\begin{gather*}
\lim_{t\rightarrow 0}\Vert \eta _{2}^{k}(t,x,\nu ,\dots ,\nu
)\Vert \leq \lim_{t\rightarrow 0}\Vert D\psi _{t}^{2}(x)\nu \Vert
=\left\Vert \nu \right\Vert
\end{gather*}
hence there is a constant $M_{k}^{\prime }>0$ such that
\begin{gather*}
\Vert D^{k}\psi _0^{2}(x)\Vert \leq M_{k}^{\prime }.
\end{gather*}

In the same way as above the integral $\int_{0}^{t}\Vert D\psi
_{t-s}^{2}(\psi _{s}^{2}(x))\Vert \Vert G_{2}^{k}(s,x,\nu )\Vert
ds$ is well def\/ined and putting $\tau =\frac{s}{t}$ we obtain
\begin{gather*}
\eta _{2}^{k}(t,x,\nu ,\dots ,\nu )=D\psi _{t}^{2}(x)\nu
+t\int_{0}^{1}D\psi _{t(1-\tau )}^{2}(\psi _{t\tau
}^{2}(x))G_{2}^{k}(t\tau ,x,\nu )d\tau .
\end{gather*}
Since $b_{1}=b^{\prime }(1+m_0)-c_{1}>a_0m_0^{\prime }$ , by the
estimates~(\ref{26}) and (\ref{28}), we deduce the existence of a constant $%
M_{k}>0$ such that
\begin{gather*}
\Vert \eta _{2}^{k}(t,x,\nu ,\dots ,\nu )\Vert \leq \Vert D\psi
_{t}^{2}(x)\nu \Vert +t\int_{0}^{1}\Vert G_{2}^{k}(t\tau ,x,\nu
)\Vert d\tau
\\
\phantom{\Vert \eta _{2}^{k}(t,x,\nu ,\dots ,\nu )\Vert}{} \leq \Vert D\psi _{t}^{2}(x)\nu \Vert +t\sum_{l=2}^{k}\int_{0}^{1}%
\frac{(t\tau )^{-\frac{lb_{1}}{a_0m_0^{\prime }}}\left\Vert
x\right\Vert ^{1+m_0^{\prime }-l}}{\left( 1+b_0m_0^{\prime }t\tau
\left\Vert
x\right\Vert ^{m_0^{\prime }}\right) ^{\frac{1+m_0^{\prime }-l}{%
m_0^{\prime }}}}d\tau \leq M_{k}t^{-\frac{b_{1}}{a_0m_0^{\prime }}}%
\text{.}
\end{gather*}
Which shows that the origin $0$ is a $G.A.S.$ equilibrium for
$\eta _{2}^{k}$ on ${\mathbb R}^{n}$. We formulate this fact as

\begin{proposition}
Let $k\geq 0$ be any integer. Under the above conditions $(i)$,
$(ii)$, $(iii)$, $(iv)$ and~$(v)$, the origin $0$ is a $G.A.S.$ of
order $k$ on ${\mathbb R}^{n}$ for the $ Y_{2}$-flow and there is
a constant $M_{k}>0$ such that for any $t\geq t_0>0$
\begin{gather}
\Vert D^{k}\psi _{t}^{2}(x)\Vert \leq M_{k}t^{-\frac{b_{1}}{%
a_0m_0^{\prime }}}. \label{33}
\end{gather}
\end{proposition}

\subsection[Global stability of prolongations of the $Y_{3}$-flow]{Global
stability of prolongations of the $\boldsymbol{Y_{3}}$-f\/low}

Let
\begin{gather*}
Y_{3}=\sum_{i=1}^{n}\big( \alpha _{i}x_{i}+\beta
_{i}x_{i}^{1+m_{i}}+Z_{3i}(x)\big) \frac{\partial }{\partial
x_{i}}
\end{gather*}
with

i) all the coef\/f\/icient $\alpha_{i}$ are negative with $-a\leq
\alpha _{i}\leq -b$;

ii) all the coef\/f\/icients $\beta_{i}\leq 0$ and $-a^{\prime
}\leq \beta _{i}\leq -b^{\prime }$;

iii) the exponents $m_{i}$ are even natural numbers with
$0<m_0\leq m_{i}\leq m_0^{\prime }$;

iv) For any $k=0,\dots,1+m_{i}$%
\begin{gather*}
\Vert D^{k}Z_{3i}(x)\Vert \leq c_{k}^{\prime }\left\vert
x_{i}\right\vert ^{2-k+m_{i}}\quad \text{if} \ x\in B\left( 0,1\right), \\
\Vert D^{k}Z_{3i}(x)\Vert \leq c_{k}^{\prime \prime }\left\vert
x_{i}\right\vert ^{1-k+m_{i}}\quad \text{if} \ x\in {\mathbb
R}^{n}\setminus B\left( 0,1\right);
\end{gather*}

v) for any $k\geq 2+m_{i}$
\begin{gather*}
\Vert D^{k}Z_{3i}(x)\Vert \leq c_{k};
\end{gather*}

vi)
\begin{gather*}
a_0=a^{\prime }+c_0, \qquad
a_{1}=a^{\prime }(1+m_0)+c_{1},\\
b_0=b^{\prime }-c_0>0, \qquad b_{1}=b^{\prime }(1+m_0)-c_{1}>0
\end{gather*}
 with $c_{k}=\max \left\{ c_{k}^{\prime },c_{k}^{\prime \prime
}\right\}.$

\begin{remark}
If $x\in B\left( 0,1\right) $ then $\Vert D^{k}Z_{3i}(x)\Vert \leq
c_{k}^{\prime }\left\vert x_{i}\right\vert ^{2-k+m_{i}}\leq
c_{k}^{\prime }\left\vert x_{i}\right\vert ^{1-k+m_{i}}$.

Let $c_{l}=\max \left\{ c_{l}^{\prime },c_{l}^{\prime \prime
}\right\} $, for any $x\in {\mathbb R}^{n}$ one has $\Vert
D^{k}Z_{3i}(x)\Vert \leq c_{k}\left\vert x_{i}\right\vert
^{1-k+m_{i}}$.
\end{remark}

\subsubsection[Global stability of the $Y_{3}$-flow ${\mathbb R}^{n}$]{Global
stability of the $\boldsymbol{Y_{3}}$-f\/low $\boldsymbol{{\mathbb
R}^{n}}$}

Denote by $\psi _{t}^{3}=\exp (tY_{3})$, by the estimates
(\ref{30}), we have
\begin{gather*}
\Vert \psi _{t}^{3}(x)\Vert \leq C \Vert x\Vert e^{-bt}\quad
\forall \; t>0 \ \text{and} \ \forall \; x\in {\mathbb R}^{n},
\end{gather*}
where $C>0$ is a constant. So by Proposition~\ref{prop1}, $0$ is a
$G.A.S.$ on ${\mathbb R}^{n}$. We proceed by recurrence; since the
property is true in case $k=0$, we assume that the property
remains true for any $l=1,\dots,k-1$, with $k$ f\/ixed i.e.\ $0$
is a global $G.A.S.$ of $\eta _{3}^{l}(t,x,\nu ,\dots \nu )=\Vert
D^{l}\psi _{t}^{3}(x)\nu ^{k}\Vert $ on ${\mathbb R}^{n}$ and
there exist constants $M_{l}>0$ such that for any $t>0$
\begin{gather*}
\Vert D^{l}\psi _{t}^{3}(x)\Vert \leq M_{l}e^{-bt}.
\end{gather*}
We will show that $0$ is a $G.A.S.$ equilibrium to $\eta _{3}^{k}$
on ${\mathbb R}^{n}$.

$\eta _{3}^{k}(t,x,\nu ,\dots,\nu )$ is a solution to the dynamic
system
\begin{gather*}
\frac{d}{dt}\eta _{3}^{k}=D_{y}\eta _{3}^{k}+G_{3}^{k}(t,x,\nu )
\end{gather*}
with $y=\psi _{t}^{3}(x)$ and
\begin{gather*}
G_{3}^{k}(t,x,\nu )=\sum_{l=2}^{k}D_{y}^{l}Y_{3}(y)\sum_{\underset{i_{j}>0}{%
i_{1}+\dots +i_{l}=k}}\left( \prod_{j=1}^{l}D^{i_{j}}\psi
_{t}^{3}(x)\nu ^{i_{j}}\right) \text{.}
\end{gather*}
By the method of the resolvent, we get
\begin{gather*}
\eta _{3}^{k}(t,x,\nu ,\dots ,\nu )=D\psi _{t}^{3}(x)\nu
+\int_{0}^{t}D\psi _{t-s}^{3}(\psi _{s}^{3}(x))G_{3}^{k}(s,x,\nu
)ds
\end{gather*}
and by the same argument as for the $Y^{1}$-f\/low, we deduce that
for any integer $k\geq $ $0$ there exist a constant $M_{k}$ such
that $\forall\; t\geq 0$
\begin{gather*}
\Vert D^{k}\psi _{t}^{3}(x)\Vert \leq M_{k}\left\Vert x\right\Vert
e^{-bt}.
\end{gather*}
By Proposition~\ref{prop1}, we have

\begin{proposition}
Under the above conditions $(i)$, $(ii)$, $(iii)$, $(iv)$, $(v)$
and $(vi)$, the origin $0$ is a $G.A.S.$ equilibrium of order $k$
on ${\mathbb R}^{n}$ to the $Y_{3}$-flow .
\end{proposition}

\pdfbookmark[1]{References}{ref}
\LastPageEnding

\end{document}